\documentclass{amsart}

\usepackage[utf8]{inputenc}
\usepackage[margin=1.2in]{geometry}
\usepackage{amsmath, amssymb, amsthm}
\usepackage{graphicx}
\usepackage{dsfont}
\usepackage{xcolor}
\usepackage{dirtytalk}
\usepackage{comment}
\usepackage{hyperref}
\usepackage{mathtools}
\usepackage{enumitem}
\usepackage[normalem]{ulem}
\usepackage{tikz}

\makeatletter
\newenvironment{customproof}[1]{%
  \par\pushQED{\qed}%
  \normalfont\topsep6pt \trivlist
  \item[\hskip\labelsep\itshape
    Proof of #1\@addpunct{.}]\ignorespaces
}{%
  \popQED\endtrivlist\@endpefalse
}
\makeatother

\theoremstyle{plain}
\newtheorem{theorem}{Theorem}[section]
\newtheorem{lemma}[theorem]{Lemma}

\newtheorem{preremark}[theorem]{Remark}

\theoremstyle{definition}
\newtheorem{definition}[theorem]{Definition}

\numberwithin{equation}{section}

\newcommand{\trace}{\operatorname{tr}}

\newcommand{\R}{\mathbb R}

\newcommand{\p}{\partial}

\newcommand{\Mp}{\mathcal{M}^+}
\newcommand{\Mm}{\mathcal{M}^-}

\title{Splitting of viscosity solutions}
\author{David Jesus}

\newcommand{\Addresses}{{
  \bigskip
  \footnotesize
    David Jesus, \textsc{Dipartimento SMFI, Universit\`a di Parma, Parco Area delle Scienze 53/a, Campus, 43124 Parma, Italy}\par\nopagebreak
  \textit{E-mail address:} {\tt davidjoao.brandligtdejesus@unipr.it}\\

 }

}

\begin{document}

\begin{abstract}
    We give a penalization argument for splitting a continuous viscosity solution of two-sided Pucci inequalities into a solution of the same inequalities with zero boundary values and a function satisfying the homogeneous Pucci inequalities and with the same boundary data, settling an open question raised by Silvestre-Sirakov.
\end{abstract}

\maketitle

\section{Introduction}

Let $\Omega\subset\R^d$ be a bounded $C^2$ domain, where $d\geq1$, and fix
$0<\lambda\leq\Lambda$ and $K\geq0$. For a symmetric matrix $M$ with eigenvalues
$\mu_1,\ldots,\mu_d$, the Pucci operators are given by
\[
 \Mp(M)=\Lambda\sum_{\mu_i>0}\mu_i
+\lambda\sum_{\mu_i<0}\mu_i,
 \qquad
 \Mm(M)=\lambda\sum_{\mu_i>0}\mu_i
+\Lambda\sum_{\mu_i<0}\mu_i.
\]
We omit the ellipticity constants from the notation. All differential
inequalities below are understood in the viscosity sense.

For $L\geq0$, we denote by $S^*_L(\Omega)$ the set of functions
$z\in C(\overline\Omega)$ satisfying
\begin{equation}\label{e:sstar}
 \mathcal M^+(D^2z)+K|Dz|\geq-L,
 \qquad
 \mathcal M^-(D^2z)-K|Dz|\leq L
 \quad\text{in }\Omega.
\end{equation}

We also consider the class with a prescribed right-hand side. For
$f\in C(\overline\Omega)$, write
\begin{equation*}
 S_f(\Omega)=\left\{z\in C(\overline\Omega):
 \mathcal M^-(D^2z)\leq f\leq\mathcal M^+(D^2z)
 \ \text{in }\Omega\right\}.
\end{equation*}
This notation does not include a gradient term. The constants $\lambda,\Lambda$ and $K$ are fixed throughout the whole paper.

To motivate our result, note that for linear elliptic equations, it is possible to separate the effects of the source term and the boundary data by subtracting the solution of the same equation with zero boundary data, and obtain a solution of the homogeneous equation. More generally, if $u$ is a viscosity solution to the fully nonlinear equation
\begin{align*}
    \begin{cases}
        F(D^2u)=f & \text{ in } \Omega,\\
        u=g & \text{ on } \p\Omega,
    \end{cases}
\end{align*}
where $F$ is uniformly elliptic with $F(0)=0$, then one can consider the solution to the same equation with zero boundary data
\begin{align*}
    \begin{cases}
        F(D^2v)=f & \text{ in } \Omega,\\
        v=0 & \text{ on } \p\Omega,
    \end{cases}
\end{align*}
and an adaptation of the argument from \cite[Theorem 5.3]{CC} gives $u-v\in S_0(\Omega)$. However, when $u$ belongs to the class $S_f$ or $S^*_L$, this property is not directly available, since it is not clear how to construct $v$. This question, concerning splitting of viscosity solutions, was  raised by Silvestre and Sirakov \cite[Section 1.1]{SS}, which we restate below.

\medskip

\noindent \textit{Open problem: } Let $u$ be a solution to $S^*$. Is it true that $u = v + w$, where $v$ solves $S^*$ and $v = 0$ on $\p\Omega$ and $w$ is a solution to $S^*$ with $L = 0$? More simply, say $u$ is a solution
of $\Mp(D^2u)\geq f(x)\geq \Mm(D^2u)$ in $\Omega$, is it true that we have the splitting $u = v + w$ where
$v$ satisfies the same inequalities and vanishes on $\p\Omega$, while $\Mp(D^2w)\geq 0\geq \Mm(D^2w)$ in $\Omega$?

\medskip

The question arose in boundary regularity theory, where such a
splitting would allow one to separate the inhomogeneous inequalities
from the nonzero boundary data, since it would considerably simplify the proof.

We answer both questions affirmatively. The first result concerns the statement in the larger class $S^*_L$.

\begin{theorem}\label{t:splitting}
Let $\Omega\subset\R^d$ be a bounded $C^2$ domain, $0<\lambda\leq\Lambda$,
and $K,L\geq0$. Suppose that $u\in S^*_L(\Omega)$. Then there are
$v,w\in C(\overline\Omega)$ such that
\begin{equation*}
 u=v+w,\qquad
 v\in S^*_L(\Omega),\qquad v=0\ \text{on }\partial\Omega,
 \qquad w\in S^*_0(\Omega),
\end{equation*}
with the same constants $\lambda,\Lambda$ and $K$.
\end{theorem}

The second answer concerns functions in $S_f(\Omega)$.

\begin{theorem}\label{t:prescribed}
Let $\Omega\subset\R^d$ be a bounded $C^2$ domain, $0<\lambda\leq\Lambda$,
and $f\in C(\overline\Omega)$. Suppose that $u\in S_f(\Omega)$. Then there are
$v,w\in C(\overline\Omega)$ such that
\begin{equation}\label{e:prescribed-splitting}
 u=v+w,\qquad
 v\in S_f(\Omega),\qquad v=0\ \text{on }\partial\Omega,
 \qquad w\in S_0(\Omega),
\end{equation}
with the same ellipticity constants.
\end{theorem}

We now discuss the main difficulties in obtaining the splittings in Theorems \ref{t:splitting} and \ref{t:prescribed}. We need to construct approximations that simultaneously 
\begin{itemize}
    \item preserve both Pucci inequalities with the same constants;
    \item attain exactly the boundary values of $u$;
    \item solve equations whose structure allows for comparison and difference arguments (that is, the difference of solutions belongs to the appropriate $S$ class).
\end{itemize}

To overcome these difficulties, we construct a new approximation method, designed to meet these three requirements together. That is the purpose of the operator $H_\varepsilon$ in Lemma \ref{l:approximation}, which is the main tool used in this paper.

\section{Preliminaries}

We denote by $ S(d)$ the space of real symmetric $d\times d$
matrices. We will often use the notation
\[
H[u](\cdot):=H(\cdot, u(\cdot),Du(\cdot),D^2u(\cdot)).
\]

\begin{definition}[Viscosity solutions]
Let $D\subset\R^d$ be open, and let
$H\in C(D\times\R\times\R^d\times S(d))$ be nondecreasing
in its matrix variable. We say $z\in USC(D)$ is a
viscosity subsolution of
\[
 H(x,z,Dz,D^2z)=0\quad\text{in }D
\]
if, whenever $x_0\in D$ and  $\phi\in C^2(B_r(x_0))$, with $B_r(x_0)\subset D$, touches  $z$ from above at $x_0$, we have
\[
 H\big(x_0,z(x_0),D\phi(x_0),D^2\phi(x_0)\big)\geq0.
\]
Similarly, $z\in LSC(D)$  is a viscosity supersolution if, whenever $x_0\in D$ and  $\phi\in C^2(B_r(x_0))$, with $B_r(x_0)\subset D$, touches  $z$ from below at $x_0$, we have
\[
 H\big(x_0,z(x_0),D\phi(x_0),D^2\phi(x_0)\big)\leq0.
\]
A function $z\in C(D)$ is a viscosity solution if it is both a
viscosity subsolution and a viscosity supersolution.
\end{definition}

We start by proving comparison, existence and uniqueness for the class of operators used later on. Since these operators are not strictly monotone in the solution value, we include some details of the proofs, using arguments from \cite{CIL}.

\begin{lemma}\label{l:dirichlet}
Let
$H\in C(\overline\Omega\times\R\times\R^d\times S(d))$
be nonincreasing in its second variable. Assume that
\begin{equation}\label{a:elliptic}
 \begin{split}
 \mathcal M^-(M-N)-K|p-q|
 &\leq H(x,r,p,M)-H(x,r,q,N)\\
 &\leq\mathcal M^+(M-N)+K|p-q|
 \end{split}
\end{equation}
for all $x,r,p,q,M,N$, and that there is a modulus of continuity $\omega$ such
that
\begin{equation}\label{a:modulus}
 |H(x,r,p,M)-H(y,r,p,M)|\leq\omega(|x-y|)
\end{equation}
uniformly in $r,p,M$. Then the comparison principle holds on every bounded open set
$D\subset\Omega$, that is, if $u\in USC(D)$ and $v\in LSC(D)$ are respectively sub and supersolutions with $u\leq v$ on $\p D$, then $u\leq v$ in $\overline D$.

If, in addition, for some $L\geq0$,
\begin{equation}\label{a:zero_bound}
 |H(x,r,0,0)|\leq L
 \quad\text{for all }x\in\overline\Omega, r\in\R,
\end{equation}
then, for every $g\in C(\partial\Omega)$, there is a unique
$z\in C(\overline\Omega)$ solving
\begin{equation}\label{e:dirichlet}
 H[z]=0\quad\text{in }\Omega,
 \qquad z=g\quad\text{on }\partial\Omega.
\end{equation}
Moreover, $z\in S^*_L(\Omega)$.
\end{lemma}

\begin{proof}
    First we prove the comparison principle. Let $u\in USC(\overline D)$ and $v\in LSC(\overline D)$ be viscosity sub/super solutions with $u\leq v$ on $\p D$. Choose $b>K/\lambda$ and $E>\max_{\overline \Omega} e^{b x_1} $, and set
    \[
    \psi(x)=e^{bx_1} -E<0.
    \]
    Since
    \begin{align*}
        \Mm(D^2\psi)-K|D\psi|=e^{bx_1}(\lambda b^2-Kb)\geq c>0\quad \text{ on } \overline \Omega,
    \end{align*}
    the function $u_\delta=u+\delta\psi$, where $\delta>0$, is a strict
subsolution satisfying
\[
 H[u_\delta]\geq\delta c\quad\text{in }D.
\]
Also, $u_\delta\leq u\leq v$ on $\partial D$.

Suppose, by contradiction, that
\[
m_\delta:=\max_{\overline D} (u_\delta-v)>0.
\]
Because of the ordering, this maximum is attained in the interior. We now double the variables in order to obtain a contradiction. For $\tau>0$, let $(x_\tau,y_\tau)$ maximize
\[
\Phi_\tau(x,y):=u_\delta(x)-v(y)-\frac{|x-y|^2}{2\tau}\quad \text{ on } \overline D\times \overline D.
\]
Clearly $\Phi_\tau(x_\tau,y_\tau)\geq m_\delta$ and
\[
u_\delta(x_\tau)-v(y_\tau)\geq m_\delta>0.
\]
Boundedness of the two functions gives
\[
|x_\tau-y_\tau|^2\leq C_\delta \tau.
\]
Therefore, every convergent subsequence of $(x_\tau,y_\tau)$, as $\tau\to 0$, has a limit of the form $(\bar x,\bar x)$. Using the respective semi-continuity of $u_\delta$ and $v$, we get
\begin{align*}
    m_\delta\leq \limsup_{\tau\to 0}(u_\delta(x_\tau)-v(y_\tau))\leq u_\delta(\bar x)-v(\bar x).
\end{align*}
Thus $\bar x\in D$ again by ordering. This implies that, for $\tau$ small enough, $x_\tau,y_\tau\in D$. At these points, the theorem of sums \cite[Theorem 3.2]{CIL} gives super and sub-jets
\begin{align*}
&(p_\tau,X_\tau)\in \overline J^{2,+}u_\delta(x_\tau),\\
   &(p_\tau,Y_\tau)\in \overline J^{2,-}v(y_\tau),
\end{align*}
with
\[
p_\tau=\frac{x_\tau-y_\tau}{\tau}.
\]
The matrices can be chosen so that
\[
 \begin{pmatrix}X_\tau&0\\0&-Y_\tau\end{pmatrix}
 \leq\frac{3}{\tau}\begin{pmatrix}I&-I\\-I&I\end{pmatrix}.
\]
Evaluating this matrix inequality on vectors of the form $(\xi,\xi)\in \R^{2d}$ shows $X_\tau\leq Y_\tau$. Therefore, since $H$ is continuous, we obtain 
\begin{align*}
    &H(x_\tau,u_\delta(x_\tau),p_\tau, X_\tau)\geq \delta c,\\
    &H(y_\tau,v(y_\tau),p_\tau, Y_\tau)\leq 0.
\end{align*}
Using $u_\delta(x_\tau)>v(y_\tau)$, monotonicity of $H$ in the function value, ellipticity and \eqref{a:modulus}, we obtain
\begin{align*}
    \delta c\,&\leq H(x_\tau,u_\delta(x_\tau),p_\tau, X_\tau)\\
    &\leq H(y_\tau,v(y_\tau),p_\tau, Y_\tau)+\omega(|x_\tau-y_\tau|)\\
    &\leq \omega(|x_\tau-y_\tau|).
\end{align*}
Letting $\tau\to 0$ contradicts $\delta c>0$. We have proved $u_\delta\leq v$ in $D$, for all $\delta>0$. Since $u_\delta\to u$ uniformly in $D$, letting $\delta\to 0$, we obtain $u\leq v$, as intended.

\medskip

We next construct boundary barriers. Fix $x_0\in \p\Omega$. Since the domain is $C^2$, there is an exterior tangent ball $B_R(y)$ whose closure meets $\overline \Omega$ only at $x_0$. Put
\[
q(x)=R^{-m}-|x-y|^{-m},\quad R_1=\max_{x\in \overline \Omega}|x-y|,
\]
where $m>0$ is chosen so that $\lambda(m+1)>\Lambda(d-1)+KR_1$. Then $q\geq 0$ on $\overline \Omega$ with equality only at $x_0$. Writing $r=|x-y|$, a direct calculation gives
\begin{equation}\label{eq:Mpq}
    \begin{aligned}
    \Mp(D^2q)+K|Dq|\,&=mr^{-m-2}\left(\Lambda(d-1)-\lambda(m+1)+Kr\right)\\
    &\leq -c_0<0.
\end{aligned}
\end{equation}
For $\eta>0$, consider
\begin{align*}
    b^+(x)=g(x_0)+\eta+Cq(x),\quad
    b^-(x)=g(x_0)-\eta-Cq(x).
\end{align*}
Since $g$ is continuous and $q$ is strictly positive away from $x_0$, we can choose $C$ so that $b^-\leq g\leq b^+$ on $\p\Omega$. Increasing $C$ to ensure that $Cc_0\geq L$, we obtain from \eqref{a:zero_bound} and \eqref{eq:Mpq} that $b^+$ and $b^-$ are super and subsolutions of \eqref{e:dirichlet}, respectively.

The comparison principle and the barriers above give existence of a
viscosity solution $z\in C(\overline\Omega)$ of \eqref{e:dirichlet} by the
Perron argument of \cite[Theorem~4.1]{CIL}, with the continuous boundary
trace $g$ ensured by these barriers. Comparison also gives uniqueness.

Finally, using properties \eqref{a:elliptic} and \eqref{a:zero_bound}, we immediately get $z\in S^*_L(\Omega)$.

\end{proof}

The following compactness statement follows from the boundary barriers and interior H\"older estimates.

\begin{lemma}\label{l:compactness}
For fixed $L\geq0$ and $g\in C(\partial\Omega)$, the family
\[
 \{z\in S^*_L(\Omega):z=g\text{ on }\partial\Omega\}
\]
is relatively compact in $C(\overline\Omega)$.
\end{lemma}

\begin{proof}
The barriers constructed in Lemma~\ref{l:dirichlet}, applied to the Pucci
inequalities for both $z$ and $-z$, give a uniform $L^\infty$ bound and a
common modulus of continuity at the boundary, depending on $g,L$ and the
fixed structural and domain parameters. Interior H\"older estimates for
$S^*_L$, obtained from the Harnack inequality in
\cite[Proposition~2.2]{SS}, give equicontinuity on compact subsets of
$\Omega$. Combining the interior and boundary estimates gives
equicontinuity on $\overline\Omega$. The conclusion follows from the
Arzel\`a--Ascoli theorem.
\end{proof}

We prove the difference property in the form needed below. For related results see  \cite[Theorem 5.3]{CC} and \cite[Theorem 2.6]{Feldman}. The assumption \eqref{a:modulus} is crucial to apply the comparison principle.

\begin{lemma}\label{l:difference}
Let $F\in C(\overline\Omega\times\R^d\times S(d))$ satisfy
\eqref{a:elliptic} and \eqref{a:modulus}, with the variable $r$
omitted. Suppose that $u,v\in C(\overline\Omega)$ solve
\[
 F(x,Du,D^2u)=F(x,Dv,D^2v)=0\quad\text{in }\Omega.
\]
Then $u-v\in S^*_0(\Omega)$.
\end{lemma}

\begin{proof}
    Let $\phi $ touch $u-v$ from above at $x_0\in \Omega$. Suppose, for the sake of contradiction, that
\[
\Mp(D^2\phi(x_0))+K|D\phi(x_0)|<0.
\]
Up to adding $|x-x_0|^4$ to $\phi$, we may suppose the contact is strict without changing the inequality above. By $\phi\in C^2$ and continuity of $\Mp$, we have for sufficiently small $r$,
\[
\Mp(D^2\phi)+K|D\phi|\leq -c<0 \quad \text{ in } B_r(x_0).
\]
We now check that $v+\phi$ is a strict supersolution of $F=0$ in $B_r(x_0)$. Indeed, if $\xi$ touches $v+\phi$ from below at $x$ in $B_r(x_0)$, then $\xi-\phi\in C^2$ touches $v$ from below and \eqref{a:elliptic} gives
\begin{align*}
    F[\xi](x)&\leq F[\xi-\phi](x)+\Mp(D^2\phi(x))+K|D\phi(x)|\\
    &\leq -c.
\end{align*}
Since the contact between $\phi$ and $u-v$ is strict, we can take $\delta>0$ such that $v+\phi-\delta\geq u$ on $\p B_r(x_0)$ while $v(x_0)+\phi(x_0)-\delta<u(x_0)$. Applying the comparison principle to the supersolution $v+\phi-\delta$ and to the subsolution $u$ in $B_r(x_0)$, we obtain a contradiction. Thus, we conclude that 
\[
\Mp(D^2(u-v))+K|D(u-v)|\geq 0.
\]
Applying the same argument to $v-u$ gives the other inequality.
\end{proof}

\section{Approximation by penalized equations}

For $p\in\R^d$ and $M\in S(d)$, set
\begin{equation*}
 A(p,M)=\mathcal M^-(M)-K|p|-L,
 \qquad
 B(p,M)=\mathcal M^+(M)+K|p|+L.
\end{equation*}
Then $A\leq B$, and \eqref{e:sstar} is equivalent to the viscosity inequalities
\[
A(Du,D^2u)\leq0\leq B(Du,D^2u).
\]

The following approximation is the main construction in the proof.
The idea behind its construction is the following. Let $u_0$ be a limit of $u_\varepsilon$. On the set $\{u>u_0\}$, the penalization term tends to $+\infty$. This selects the branch $B$ in the viscosity argument, giving $B[u_0]\leq0$ in this set. Since $B[u]\geq0$ and $u=u_0$ on its boundary, the comparison principle implies $u\leq u_0$, a contradiction.
On the set $\{u<u_0\}$, the penalization instead tends to $-\infty$, and the same argument using the branch $A$ gives the opposite inequality. Thus $u_\varepsilon\longrightarrow  u$. Moreover, as we will see, the operator $H_\varepsilon$  allows for splitting of $u_\varepsilon$.

\begin{lemma}\label{l:approximation}
Let $u\in S^*_L(\Omega)$ and fix $a\in[\lambda,\Lambda]$. For $\varepsilon>0$, define
\begin{equation*}
 H_\varepsilon(x,r,p,M)
 =\max\left\{A(p,M),
\min\left\{B(p,M),a\trace M+
\frac{u(x)-r}{\varepsilon}\right\}\right\}.
\end{equation*}
There is a unique $u_\varepsilon\in C(\overline\Omega)$ solving
\begin{equation}\label{e:penalized-problem}
 H_\varepsilon[u_\varepsilon]=0\quad\text{in }\Omega,
 \qquad u_\varepsilon=u\quad\text{on }\partial\Omega.
\end{equation}
Moreover, $u_\varepsilon\in S^*_L(\Omega)$ and
\begin{equation}\label{e:uniform-approximation}
 u_\varepsilon\longrightarrow u\quad\text{uniformly on }\overline\Omega
 \quad\text{as }\varepsilon\to 0.
\end{equation}
\end{lemma}

\begin{proof}
It is simple to check that $H_\varepsilon$ satisfies \eqref{a:elliptic} since each of its three defining branches also does, and taking a maximum or minimum preserves this property. Moreover, it is continuous and nonincreasing in $r$.

Let $\omega_u$ be a modulus of $u$ on $\overline \Omega$. Then
\begin{align*}
    |H_\varepsilon(x,r,p,M)-H_\varepsilon(y,r,p,M)|\leq \frac{\omega_u(|x-y|)}{\varepsilon}.
\end{align*}
Also, $|H_\varepsilon(x,r,0,0)|\leq L$. Lemma \ref{l:dirichlet} therefore gives existence and uniqueness of a solution $u_\varepsilon$ to \eqref{e:penalized-problem}. Since $A\leq H_\varepsilon\leq B$, we have $u_\varepsilon\in S^*_L(\Omega)$. 

By Lemma \ref{l:compactness}, every sequence $\varepsilon_j\to 0$ has a subsequence, which we denote again by $\varepsilon_j$, such that $u_{\varepsilon_j}$ converges uniformly on $\overline \Omega$ to a function $u_0$ coinciding with $u$ on $\p\Omega$. We prove $u_0=u$ in the interior as well.

Consider the open set 
\[
D_+=\{x\in \Omega\,:\, u_0(x)>u(x)\}.
\]
On every compact subset of $D_+$, we have
\begin{equation}\label{eq:blow-up}
    \frac{u-u_{\varepsilon_j}}{\varepsilon_j}\longrightarrow -\infty\quad \text{ uniformly }.
\end{equation}
Let $x_0\in D_+$ and $\phi\in C^2$ touch $u_0$ from above at $x_0$ in $D_+$. We can suppose the contact is strict, up to adding $|x-x_0|^4$. By uniform convergence of $u_{\varepsilon_j}$ to $u_0$, there are points $x_j$ converging to $x_0$ such that $u_{\varepsilon_j}-\phi$ has a local maximum at $x_j$. Moreover, for $j$ large enough, all these points are contained in $D_+$, so \eqref{eq:blow-up} holds at $x_j$. Since $u_{\varepsilon_j}$ solves \eqref{e:penalized-problem} with $\varepsilon=\varepsilon_j$, we get
\begin{align*}
    0\,&\leq H_{\varepsilon_j}(x_j, u_{\varepsilon_j}(x_j),D\phi(x_j),D^2\phi(x_j)) \\
    &=\max\left\{  A[\phi](x_j),\min\left\{ B[\phi](x_j),a \trace D^2\phi(x_j)+\frac{u(x_j)-u_{\varepsilon_j}(x_j)}{\varepsilon_j}  \right\}  \right\}.
\end{align*}
Since $|D^2\phi(x_j)|, |D\phi(x_j)|\leq C$ for some $C>0$, then for $j$ large enough
\[
H_{\varepsilon_j}(x_j, u_{\varepsilon_j}(x_j),D\phi(x_j),D^2\phi(x_j))=A[\phi](x_j)\geq 0.
\]
Passing to the limit, we get
\[
A[\phi](x_0)\geq 0
\]
which implies
\[
A[u_0]\geq 0 \quad \text{ in } D_+.
\]
On the other hand, since $u\in S^*_L(\Omega)$, then $u$ is a supersolution of this equation. Continuity implies $u_0=u$ on $\p D_+$. Comparison on $D_+$ therefore gives $u_0\leq u$ so $D_+$ is empty.

On $D_-=\{x\in \Omega\,:\, u_0(x)<u(x)\}$, the penalty instead tends to $+\infty$, hence repeating the same argument with obvious changes gives
\[
B[u_0]\leq 0 \quad \text{ in } D_-,
\]
and comparison with the subsolution $u$, for which $B[u]\geq 0$, gives $u\leq u_0$ in $D_-$. Thus $D_-$ is empty.

We proved that every sequence $u_{\varepsilon_j}$, with $\varepsilon_j\to 0$, has a subsequence converging uniformly to $u$ on $\overline \Omega$. This implies the convergence of the family $u_\varepsilon$.

\end{proof}

\section{Proof of the splitting}

\begin{customproof}{Theorem \ref{t:splitting}}
    Fix $a\in [\lambda,\Lambda]$, independently of $\varepsilon$, and let $u_\varepsilon$ be the corresponding solution from Lemma \ref{l:approximation}. Define
    \[
b_\varepsilon(x)=\frac{u(x)-u_\varepsilon(x)}{\varepsilon}
    \]
and freeze the coefficient, setting
\begin{equation*}
    F_\varepsilon(x,p,M)=\max\left\{ A(p,M), \min\left\{ B(p,M), a\trace M+b_\varepsilon(x)   \right\}    \right\}.
\end{equation*}
Then $u_\varepsilon$ satisfies, in the viscosity sense,
\begin{equation*}
    F_\varepsilon(x,Du_\varepsilon,D^2u_\varepsilon)=0 \quad \text{ in } \Omega.
\end{equation*}
Indeed, if $\phi$ touches $u_\varepsilon$ at $x_0$ from either side, then $\phi(x_0)=u_\varepsilon(x_0)$ and we obtain exactly $b_\varepsilon(x_0)$ in the third branch.

The operator $F_\varepsilon$ satisfies \eqref{a:elliptic} and
\begin{align*}
|F_\varepsilon(x,p,M)-F_\varepsilon(y,p,M)|\leq |b_\varepsilon(x)-b_\varepsilon(y)|,\quad |F_\varepsilon(x,0,0)|\leq L.
\end{align*}

For each fixed $\varepsilon$, uniform continuity of $b_\varepsilon$ gives a modulus of continuity independent of $p,M$. By Lemma \ref{l:dirichlet}, there is a solution $v_\varepsilon\in C(\overline \Omega)$ of
\begin{align*}
    F_\varepsilon(x,Dv_\varepsilon,D^2v_\varepsilon)=0 \quad \text{ in } \Omega,\quad v_\varepsilon=0 \quad \text{ on } \p\Omega.
\end{align*}
Since $A\leq F_\varepsilon\leq B$, we have $v_\varepsilon\in S_L^*(\Omega)$. Since both $u_\varepsilon$ and $v_\varepsilon$ solve $F_\varepsilon=0$, we can apply Lemma \ref{l:difference} to get
\[
w_\varepsilon=u_\varepsilon-v_\varepsilon\in S_0^*(\Omega).
\]
The functions $v_\varepsilon$ have zero boundary values, so Lemma \ref{l:compactness} gives a sequence $\varepsilon_j\to 0$ and $v\in C(\overline \Omega)$ such that $v_{\varepsilon_j}\to v$ uniformly on $\overline \Omega$. By \eqref{e:uniform-approximation}, 
\[
w_{\varepsilon_j} \longrightarrow w:=u-v \quad \text{ uniformly on } \overline \Omega.
\]
Viscosity inequalities with continuous operators are stable under uniform convergence. Therefore, we conclude
\[
v\in S_L^*(\Omega),\quad w\in S_0^*(\Omega),\quad v=0\quad \text{ on } \p\Omega.
\]
Since $u=v+w$, the proof is complete.

\end{customproof}

The argument for Theorem \ref{t:prescribed} is similar to the one above.

\begin{customproof}{Theorem \ref{t:prescribed}}

Fix $a\in [\lambda,\Lambda]$ and set
\[
G_\varepsilon(x,r,M)=\max\left\{  \Mm(M),\min\left\{ \Mp(M),a\trace M+\frac{u(x)-r}{\varepsilon}\right\}\right\}.
\]
In order to apply Lemma \ref{l:dirichlet} to $G_\varepsilon-f$ with $K=0$ and $L=\|f\|_{L^\infty(\Omega)}$, we check that it satisfies the assumptions. First, \eqref{a:elliptic} and \eqref{a:zero_bound} follow directly from the definition. Moreover, \eqref{a:modulus} follows with $\omega=\varepsilon^{-1}\omega_u+\omega_f$. Hence, we obtain a unique $u_\varepsilon\in C(\overline \Omega)$ solving 
\[
G_\varepsilon[u_\varepsilon]=f \quad \text{ in } \Omega,\quad  u_\varepsilon=u
\quad \text{ on } \p\Omega.
\]
Moreover, $u_\varepsilon\in S_f(\Omega)\subset S_{\|f\|_\infty}^*(\Omega)$ with $K=0$.

The proof of Lemma \ref{l:approximation} gives $u_\varepsilon\to u$ uniformly, arguing exactly in the same way except using in the end the fact that $u\in S_f(\Omega)$. With $b_\varepsilon$ defined above, put
\[
Q_\varepsilon(x,M)=\max\left\{  \Mm(M),\min\left\{ \Mp(M),a\trace M+b_\varepsilon(x)\right\}\right\}.
\]
Then $Q_\varepsilon(x,0)=0$ and $Q_\varepsilon(x,D^2u_\varepsilon)=f$. Let $v_\varepsilon$ solve
\[
Q_\varepsilon(x,D^2v_\varepsilon)=f\quad \text{ in } \Omega,\quad v_\varepsilon=0\quad \text{ on } \p\Omega.
\]
We have $v_\varepsilon\in S_f(\Omega)$. Applying Lemma \ref{l:difference} to the operator $Q_\varepsilon-f$, gives 
\[
w_\varepsilon:=u_\varepsilon-v_\varepsilon\in S_0(\Omega).
\]

As in the proof of Theorem \ref{t:splitting}, we use compactness and stability to get the splitting in \eqref{e:prescribed-splitting}.

\end{customproof}

\bigskip

\noindent\textbf{Acknowledgements. } D.J. is supported by the European Research Council, through the ERC StG project NEW, No.~101220121.

\bigskip

\noindent\textbf{AI declaration. } 
The author used ChatGPT Astra 6.0 to assist with developing and checking mathematical arguments, including the penalized construction in Lemma \ref{l:approximation}, and with drafting and revising the manuscript. The author has verified the mathematical arguments and takes full responsibility for the contents.
\medskip

\noindent{\bf Data Availability.} Not applicable. 

\medskip
\noindent{\bf Conflict of interest.} The author states that there is no conflict of interest.

\bibliographystyle{amsalpha}

{\raggedright\Addresses\par}

\end{document}